\documentclass{bcp}
\usepackage{amssymb}

\def \one
{{\leavevmode\rlap{\rlap{$\upharpoonleft$}\kern.18em]}\phantom{$\upharpoonleft$}}}
\begin{document}

\keywords{eigenvalue, eigenfunction, integral operator}
\mathclass{Primary 47A75; Secondary 47A10, 47G10.}
\thanks{This research was partially supported by NAS of Ukraine, Grant  \# 0105U006289.}
\abbrevauthors{I.Yu.Domanov} \abbrevtitle{On the spectrum  of  the
operator
$(V\lowercase{f})(\lowercase{x})=\int_0^{\lowercase{x}^\alpha}\lowercase{f}(\lowercase{t})\lowercase{dt}$}

\title{ On the spectrum and eigenfunctions\\ of the operator
$(V\lowercase{f})(\lowercase{x})=\int_0^{\lowercase{x}^\alpha}\lowercase{f}(\lowercase{t})\lowercase{dt}$}

\author{I.Yu.Domanov}
\address{Institute of Applied Mathematics and Mechanics NAS
of Ukraine\\
Roza-Luxemburg str 74, Donetsk, 83114, Ukraine\\
 E-mail:  domanovi@yahoo.com}

\maketitlebcp


\section*{1. Introduction.} It is well known that the Volterra operator
$V:\ f\rightarrow \int_0^xf(t)dt$ defined on $L^p(0,1)$ $(C[0,1])$
is quasinilpotent, that is $\sigma(V)=\{0\}$. It was pointed out
in [5]-[6] that the operator
\begin{equation}
V_{\phi}:\ f\rightarrow\int_0^{\phi(x)}f(t)dt
\end{equation}
which is a composition of integration and  substitution  with
$\phi\in C[0,1]$   is quasinilpotent on $C[0,1]$ if $\phi(x)\leq
x$ for all $x\in [0,1]$.

Let $\phi :\  [0,1]\longrightarrow [0,1]$ be a measurable function
and $V_{\phi}:\ L^p(0,1)\longrightarrow L^p(0,1)$ $(1\leq
p<\infty)$ be  defined by (1) . It was proved in [12]-[13] that
$V_{\phi}$ is quasinilpotent on $L^p(0,1)$ if and only if
$\phi(x)\leq x$ for almost all $x\in [0,1]$. It was also noted in
[13] and proved in [14] that the spectral radius of $V_{x^\alpha}$
defined on $L^p(0,1)$ or $C[0,1]$  is $1-\alpha$ $(0<\alpha<1)$.

We note also  paper [4], where the hypercyclicity of
$V_{x^\alpha}$ was proved on some Fr\'{e}chet space.

In this note we find the  spectrum of $V_{x^\alpha}$ defined on
$L^2(0,1)$ and investigate some properties of its eigenfunctions.

 {\bf Notations:}
Let $X$ be a Banach space and let $T$ be a bounded operator on
$X$. Then   ${\rm ker} T:=\{x\in X\ :\ Tx=0\}$  denotes a kernel
of $T$ and  ${\rm R}(T):=\{Tx\ :\ x\in X\}$ denotes a range of
$T$.\  $I$ denotes the identity operator on $X$; ${\rm span} E$
denotes the closed linear span of the set $E\subset X$;\
{{\leavevmode\rlap{\rlap{$\upharpoonleft$}\kern.18em]}\phantom{$\upharpoonleft$}}}
 denotes the function $f\equiv1$ in $L^2(0,1)$;
 $ \mathbb Z_+ := \{0,1,2,\dots\}$.
 For simplicity we set
$\sum_{k=n}^ma_k:=0$ if $n>m$.

\section*{2. Auxiliary results.}
The following two Lemmas are well known. For the sake of
completeness,  proofs are given.

\th{lemma}{1.}
{%
 The system $\{(\ln x)^n \}_{n=0}^\infty$ is complete  in $L^2(0,1)$.
}

\begin{proof}
 Since
 the Laguerre functions
$f_n(x):=e^{-x/2}\frac{1}{n!}e^x\frac{d^n}{dx^n}(x^ne^{-x})$
$(n\in \mathbb Z_+ )$ form [1] an orthonormal basis in
$L^2(0,\infty)$, the
 system $\{x^ne^{-x/2}\}_{n=0}^\infty$ is complete  in
 $L^2(0,\infty)$.
Let the operator $T:\ L^2(0,\infty) \longrightarrow L^2(0,1)$ be
defined by
$$
(Tf)(x): =\frac{f(-\ln x)}{x^{1/2}}.
$$
It is easily proved that  $T$ is a surjective isometry. Thus the
system $\{T(x^ne^{-x/2})\}_{n=0}^\infty=\{(-\ln
x)^n\}_{n=0}^\infty$ is  complete  in $L^2(0,1)$.
\end{proof}

\remar{Remark\ {1.}\ }{ Consider an operator $C:\
L^2(0,1)\rightarrow L^2(0,1)$  defined by
$(Cf)(x)=f(x)-\int_x^1\frac{f(t)}{t}dt$.
  It is well known [2] that $C$  is a simple unilateral shift.
  Since ${\rm ker} C^*=\{c\cdot$
{\leavevmode\rlap{\rlap{$\upharpoonleft$}\kern.18em]}\phantom{$\upharpoonleft$}}
: $c\in\mathbb C\}$, it follows [8] that    the set $\{C^n$
{\leavevmode\rlap{\rlap{$\upharpoonleft$}\kern.18em]}\phantom{$\upharpoonleft$}}
$\}_{n=0}^\infty$   forms an orthonormal basis in $L^2(0,1)$.  It
can easily be checked that $(C^n$
{\leavevmode\rlap{\rlap{$\upharpoonleft$}\kern.18em]}\phantom{$\upharpoonleft$}}$)(x)$
$=P_n(\ln x)$, where $P_n$ is a polynomial of degree $n$.
 Thus $L^2(0,1)={\rm span}\{(\ln
x)^n:n\geqslant 0\}$.}

\th{lemma}{2} {Let $A$ be a compact operator defined on a Hilbert
space $H$, $Af_n=\lambda_nf_n$ and ${\rm span}\{f_n: n\geqslant
1\}=H$. Then

1) $\sigma_p(A)=\{\lambda_n\}_{n=1}^\infty$;

2) if $\lambda_i\ne\lambda_j$ for $i\ne j$ then
 for every eigenvalue of $A$ the algebraic multiplicity is equal
 to one.}
\begin{proof}
 {\bf 1)}Let $\lambda\in \sigma_p(A)$ and
$\lambda\ne\lambda_n$ for all $n=1,2,\dots$. Then
$\overline{\lambda}\in \sigma_p(A^*)$ and hence
$$
H\ne\left({\rm ker}(A^*-\overline{\lambda} I)\right)^\perp\
=\overline{{\rm R}(A-\lambda I)}={\rm span}\{(A-\lambda I)f_n:
n\geqslant 1\}
$$
$$
={\rm span}\{(\lambda_n-\lambda)f_n: n\geqslant 1\}={\rm
span}\{f_n: n\geqslant 1\}=H.
$$
This contradiction proves 1).

{\bf 2)} Let $\lambda_k\in\sigma_p(A)$. Since $A$ is a compact
operator and ${\rm span}\{f_n: n\geqslant 1\}=H$, we obtain
$$
{\rm dim} {\rm ker}(A-\lambda_k I)^m={\rm dim} \overline{{\rm
R}(A-\lambda_k I)^m} ^\perp={\rm dim} \left({\rm
span}\{(\lambda_n-\lambda_k)^mf_n:n\ge 0\}\right)^\perp=
$$
$$
{\rm dim} \left({\rm span}\{f_n:n\ge 0,\ n\ne
k\}\right)^\perp=1,\qquad\qquad m=1,2,\dots.
$$
Hence the algebraic multiplicity of $\lambda_k$ is equal to  one.
\end{proof}

 The following Lemma is a rephrasing of the
Problems I.50, V.161, V.162 from  [9].
 \th{lemma}{3.} {Let $|q|<1$
then

1) $
F_q(z):=\prod\limits_{k=1}^\infty(1-q^kz)=1+\sum\limits_{k=1}^\infty\frac{q^{k(k+1)/2}}{(q-1)\cdots(q^k-1)}z^k
 $
is an entire function.

 2) The polynomials
$P_n(z):=1+\sum\limits_{k=1}^n\frac{n!}{(n-k)!}\frac{q^{k(k+1)/2}}{(q-1)\cdots(q^k-1)}z^k$
have only real positive zeroes.}

\section*{3. Main results.}\qquad

\th{proposition}{1.} {Let $0<\alpha<1$ and
$V_\alpha:=V_{x^\alpha}$ be defined on $L^2(0,1)$. Then

 1) $\sigma_p(V_\alpha)=\{(1-\alpha)\alpha^{n-1}\}_{n=1}^{\infty}$;

 2) the algebraic multiplicity of  every eigenvalue of $V_\alpha$ is equal to one;

 3)
 $$
 f_{n+1}(x)=x^{\frac{\alpha}{1-\alpha}}\left(\ln^n x+ \sum\limits_{k=1}^{n}\frac{n!}{(n-k)!}
 \frac{\alpha^{k(k-1)/2}(1-\alpha)^k}{(1-\alpha)\dots
 (1-\alpha^k)}\ln^{n-k}x\right), \qquad n\in
\mathbb Z_+
 $$
 is an eigenfunction  for the operator  $V_\alpha$  with eigenvalue
 $\lambda_{n+1}:=(1-\alpha)\alpha^n$;

 4)
 $$
 g_{n+1}(x)=
 1+\sum\limits_{k=2}^\infty(-1)^{k-1}\frac{\alpha^{(k-1)(k-2-2n)/2}}{(1-\alpha)\dots
 (1-\alpha^{k-1})}x^{\frac{1-\alpha^{k-1}}{(1-\alpha)\alpha^{k-1}}}, \qquad n\in
\mathbb Z_+
 $$
  is an eigenfunction for the operator  $V_\alpha^*$  with eigenvalue
 $\lambda_{n+1}:=(1-\alpha)\alpha^n$.

 5) the system $\{f_n\}_{n=1}^\infty$ is complete in $L^2(0,1)$;

 6) the system $\{g_n\}_{n=1}^\infty$ is not complete in  $L^2(0,1)$.

 7) the operator $V_\alpha$ does not admit a spectral synthesis,
 i.e.  there exists an invariant subspace $E$ such that
 $V_\alpha|_E$ is quasinilpotent.}

\begin{proof}
 {\bf 3)} Since $x^\varepsilon\ln^m x\in C[0,1]$ for all
$\varepsilon>0$ and $m\in\mathbb Z_+$, we have that $f_{n+1}\in
L^2(0,1)$. Let us check that $f_{n+1}(x)$ is an eigenfunction of
$V_\alpha$ corresponding to the eigenvalue
 $\lambda_{n+1}:=(1-\alpha)\alpha^n$. By definition, put
$$
C_{n-k}(\alpha):=\frac{n!}{(n-k)!}
 \frac{\alpha^{k(k-1)/2}(1-\alpha)^k}{(1-\alpha)\dots
 (1-\alpha^k)},\qquad\qquad\qquad k=1\dots n.
$$
Then
$$
\frac{\alpha}{1-\alpha}C_{n-k}(\alpha)+(n-k+1)C_{n-k+1}(\alpha)=
\frac{n!}{(n-k)!}
 \frac{\alpha^{(k-1)(k-2)/2}(1-\alpha)^{k-1}}{(1-\alpha)\dots
 (1-\alpha^{k-1})}\left(\frac{\alpha^k}{1-\alpha^k}+1\right)
$$
$$
=\frac{n!}{(n-k)!}
 \frac{\alpha^{(k-1)(k-2)/2}(1-\alpha)^{k-1}}{(1-\alpha)\dots
 (1-\alpha^k)},\qquad\qquad\qquad k=1\dots n.
$$
Further,
$$
\alpha x^{\alpha-1}f_{n+1}(x^\alpha )=\alpha
x^{\alpha-1}(x^\alpha)^\frac{\alpha}{1-\alpha}\left(\ln^n
x^\alpha+
\sum\limits_{k=1}^nC_{n-k}(\alpha)\ln^{n-k}x^\alpha\right)
$$
\begin{equation}
 = \alpha x^{\alpha-1+\frac{\alpha^2}{1-\alpha}}\left(
\alpha^n\ln^n x+\sum\limits_{k=1}^n\frac{n!}{(n-k)!}
 \frac{\alpha^{k(k-1)/2}(1-\alpha)^k}{(1-\alpha)\dots
 (1-\alpha^k)}\alpha^{n-k}\ln^{n-k}x\right)
\end{equation}
$$
 =
(1-\alpha)\alpha^nx^\frac{2\alpha-1}{1-\alpha}\left(\frac{\alpha\ln^n
x}{1-\alpha}+ \sum\limits_{k=1}^n\frac{n!}{(n-k)!}
\frac{\alpha^{(k-1)(k-2)/2}(1-\alpha)^{k-1}}{(1-\alpha)\dots
(1-\alpha^k)}\ln^{n-k}x\right), \quad n\in \mathbb Z_+,
$$
and
$$
f'_{n+1}(x)=\frac{\alpha}{1-\alpha}x^{\frac{\alpha}{1-\alpha}-1}\left(\ln^n
x+ \sum\limits_{k=1}^n C_{n-k}(\alpha)\ln^{n-k}x\right)
$$
\begin{equation}
 + x^\frac{\alpha}{1-\alpha}\left(\frac{n\ln^{n-1}x}{x}+
\sum\limits_{k=1}^{n-1}
C_{n-k}(\alpha)\frac{1}{x}(n-k)\ln^{n-k-1}x\right)
\end{equation}
$$
= x^\frac{2\alpha-1}{1-\alpha}\left(\frac{\alpha\ln^n
x}{1-\alpha}+n\ln^{n-1}x\right)
$$
$$
+x^\frac{2\alpha-1}{1-\alpha} \left(\sum\limits_{k=1}^n
\frac{\alpha C_{n-k}(\alpha)}{1-\alpha}\ln^{n-k}x+
\sum\limits_{k=2}^{n} C_{n-k+1}(\alpha)(n-k+1)\ln^{n-k}x \right)
$$
$$
= x^\frac{2\alpha-1}{1-\alpha} \left[\frac{\alpha\ln^n
x}{1-\alpha}+ \frac{n}{1-\alpha}\ln^{n-1}x+\sum_{k=2}^n\left(
\frac{\alpha
C_{n-k}(\alpha)}{1-\alpha}+(n-k+1)C_{n-k+1}(\alpha)\right)\ln^{n-k}x\right]
$$
$$
= x^\frac{2\alpha-1}{1-\alpha}\left( \frac{\alpha\ln^n
x}{1-\alpha}+\sum_{k=1}^n \frac{n!}{(n-k)!}
 \frac{\alpha^{(k-1)(k-2)/2}(1-\alpha)^{k-1}}{(1-\alpha)\dots
 (1-\alpha^k)}\ln^{n-k}x\right), \qquad n\in
\mathbb Z_+.
$$
It follows from  (2)-(3)  that $ \alpha
x^{\alpha-1}f_{n+1}(x^\alpha)=(1-\alpha)\alpha^nf_{n+1}'(x)$. Thus
$$
(V_\alpha f_{n+1})(x)=
 \int_0^{x^\alpha}f_{n+1}(t)dt=
  \int_0^x\alpha t^{\alpha-1}f_{n+1}(t^\alpha)dt=
  (1-\alpha)\alpha^n\int_0^xf_{n+1}'(t)dt
$$
$$
=(1-\alpha)\alpha^n(f_{n+1}(x)-f_{n+1}(0))=
(1-\alpha)\alpha^nf_{n+1}(x), \qquad\qquad\qquad n\in \mathbb Z_+.
$$
{\bf 4)} The  convergence of the series
$$
S:=\sum\limits_{k=2}^\infty\frac{\alpha^{(k-1)(k-2-2n)/2}}{(1-\alpha)\dots
 (1-\alpha^{k-1})}x^{k-1},\qquad x\in [0,1]
$$
follows from  D'Alembert rule. Since
$\frac{\alpha^{k-1}-1}{(\alpha-1)(\alpha^{k-1})}=\frac{1}{\alpha}+\dots\frac{1}{\alpha^{k-1}}>k-1$,
we obtain that
$x^{k-1}>x^{\frac{\alpha^{k-1}-1}{(\alpha-1)(\alpha^{k-1})}}$ for
$x\in [0,1]$. Now the absolute convergence of $g_n(x)$ for $x\in
[0,1]$ ( and hence continuity of $g_n$) is implied by the
convergence of $S$.

 Let us check that $g_{n+1}(x)$ is an eigenfunction for the
 operator $V_\alpha^*$  with a corresponding eigenvalue
 $\lambda_{n+1}:=(1-\alpha)\alpha^n$.
$$
 (V_\alpha^*g_{n+1})(x)
 =\int_{x^{1/\alpha}}^1g_{n+1}(t)dt
 $$
\begin{equation}
 =1-x^{1/\alpha}+\sum\limits_{k=2}^\infty\frac{(-1)^{k-1}\alpha^{(k-1)(k-2-2n)/2}}{(1-\alpha)\dots
 (1-\alpha^{k-1})}\frac{(1-\alpha)\alpha^{k-1}}{1-\alpha^k}x^{\frac{1-\alpha^k}{(1-\alpha)\alpha^{k-1}}}\Big\vert_{x^{1/\alpha}}^1
 \end{equation}
 $$
=(1-\alpha)\alpha^n\sum\limits_{k=1}^\infty\frac{(-1)^{k-1}\alpha^{k(k-1-2n)/2}}{(1-\alpha)\dots
 (1-\alpha^k)}\left(1-x^{\frac{1-\alpha^k}{(1-\alpha)\alpha^{k-1}}}\right)=:\lambda_{n+1}(S_1-S_2)
 $$
 $$
=\lambda_{n+1}(S_1-(1-g_{n+1}(x)))=\lambda_{n+1}(S_1-1)+\lambda_{n+1}g_{n+1}(x).
 $$
By   Lemma 3 1)
$$
S_1=\sum\limits_{k=1}^\infty\frac{(-1)^{k-1}\alpha^{k(k-1-2n)/2}}{(1-\alpha)\dots
 (1-\alpha^k)}=-\sum\limits_{k=1}^\infty\frac{\alpha^{k(k+1)/2}\alpha^{(-n-1)k}}{(\alpha-1)\dots
 (\alpha^k-1)}
$$
\begin{equation}
=-(F_\alpha(\alpha^{-n-1})-1)=1.
 \end{equation}
Combining   (4) and (5), we get
$(V_\alpha^*g_{n+1})(x)=\lambda_{n+1}g_{n+1}(x)$.

{\bf 5)}   It can be proved that $ E_{n+1}:={\rm
span}\{f_1,\dots,f_{n+1}\} ={\rm
span}\{x^{\frac{\alpha}{1-\alpha}}\ln^k x:k=0\dots n\}$. Hence by
Lemma 1
$$
E_\infty:={\rm span}\{f_k:k\in \mathbb Z_+\}={\rm
span}\{x^{\frac{\alpha}{1-\alpha}}\ln^k x:k\in \mathbb
Z_+\}=\overline{x^{\frac{\alpha}{1-\alpha}}L^2(0,1)}=L^2(0,1).
$$
{\bf 1)}, {\bf 2)}  follow from 5) and Lemma 2.

{\bf 6)} It follows from  M\H{u}ntz-Sz\'{a}sz theorem [7],[11] that
the system
$\{x^{\frac{1-\alpha^n}{(1-\alpha)\alpha^n}}\}_{n=0}^\infty$ is
not complete  in $L^2(0,1)$.
 Since ${\rm span}\{g_n: n\ge 1\}\subset{\rm
span}\{x^{\frac{1-\alpha^n}{(1-\alpha)\alpha^n}}:n\ge 0\}$, we
have that  the system $\{g_n\}_{n=1}^\infty$ is not complete in
$L^2(0,1)$.

{\bf 7)} Let $E={\rm span}\{g_n: n\ge 1\}^\perp$. Then $V_\alpha
E\subset E$ and by 5) operator $V_\alpha|_E$ is quasinilpotent.
\end{proof}

\th{corollary}{1.} {Let $0<\alpha<1$,
$\phi(x)=1-(1-x)^{1/\alpha}$. Then
 operators $V_{x^\alpha}^*$ and $V_{\phi}$ are unitarily equivalent
 and hence
 $\sigma_p(V_\phi)=\{(1-\alpha)\alpha^{n-1}\}_{n=1}^{\infty}$.}

\begin{proof}
Let $U$ be a unitary operator defined by $(Uf)(x)=f(1-x)$. Then
simple computations show that $V_{x^\alpha}^*=U^{-1}V_\phi U$. 
\end{proof}

\remar{Remark\ {2.}\ }{ Suppose $\phi(x)=(1-(1-x)^{1/\alpha})'$,
then $\phi'(0)=1/\alpha$. Thus  Corollary 1 states that
 condition  $\phi'(0)=\infty$ is not necessary for ${\rm card} \{\sigma_p(V_\phi)\}=\infty$.
}

\remar{Remark\ {3.}\ } {It is interesting to note that if
$\phi(\phi(x))=x$ then the operator $V_\phi$ is selfadjoint, and
hence eigenfunctions of $V_\phi$ form an orthonormal basis in
$L^2(0,1)$. The statements 5) and 6) of Proposition 1 imply that
the operator $V_\alpha$ is not similar and even quasisimilar( see
definition in [8], [10]) to $V_\alpha^*$. It contrasts to the case
$\alpha=1$ : $V^*=U^{-1}VU$.

It follows also that $V_\alpha$ is not quasisimilar to any
selfadjoint operator.}

\th{corollary}{2.} { 1) $f_n(x)$ is a continuous function with $n$
real zeroes  which belong to $[0,1]$;

2) zeroes of $f_n(x)$ and  $f_{n+1}(x)$ interlace.
 }

\begin{proof}
 {\bf 1.} The continuity of  $f_n(x)$ was proved in
Proposition 1. Let us prove that the function $f_{n+1}$ has $n+1$
zeroes which belong to $[0,1]$. By definition, put
$$
P_n(x):=\left(\frac{t^{-\frac{\alpha}{1-\alpha}}f_{n+1}(t)}{\ln^nt}\Big\vert_{t=e^{-\frac{1-\alpha}{\alpha
x}}} \right)
$$
$$
=\left(1+\sum\limits_{k=1}^\infty\frac{n!}{(n-k)!}
\frac{\alpha^{k(k-1)/2}(1-\alpha^k)}{(1-\alpha)\dots(1-\alpha^k)}\ln^{-k}t\right)\Big\vert_{t=e^{-\frac{1-\alpha}{\alpha
x}}}
$$
$$
=1+\sum\limits_{k=1}^\infty\frac{n!}{(n-k)!}
\frac{\alpha^{k(k+1)/2}}{(\alpha-1)\dots(\alpha^k-1)}x^k.
$$
It can easily be checked that
$$
f_{n+1}(t)=t^{\frac{\alpha}{1-\alpha}}\ln^n
tP_n\left(\frac{-\alpha}{(1-\alpha)\ln t}\right).
$$
It follows from  Lemma 3 2) that the polynomial $P_n$ has exactly
$n$ positive zeroes. Thus the function $f_{n+1}$ has $n+1$ zeroes
which belong $[0,1]$.

{\bf 2.} Let us note that
$(x^nP_{n+1}(x^{-1}))'=nx^{n-1}P_n(x^{-1})$. Therefore zeroes of
$P_n(x)$ and $P_{n+1}(x)$ interlace. Hence zeroes of $f_n(x)$ and
$f_{n+1}(x)$ interlace. 
\end{proof}

\remar{Remark\ {4.}\ }{
 We suppose that eigenfunctions $g_n$ of
the operator $V_{x^\alpha}^*$ have the same  properties of zeroes
as $f_n$. Namely

1) $g_n(x)$ is a continuous function with $n$ real zeroes which
belong to $[0,1]$;

2) zeroes of $g_n(x)$ and  $g_{n+1}(x)$ interlace. }

\remar{Remark\ {5.}\ } {Proposition 1 as well as Corollary 2 hold
also if the operator $V_\alpha$ is defined on $L^p(0,1)$ ($1\le
p<\infty$). To prove it one can easily check that the operator
$V_\alpha$ defined on $L^2(0,1)$ is quasisimilar to the operator
$V_\alpha$ defined on $L^p(0,1)$.}

\remar{Remark\ {6.}\ } {It was proved in  [4] that $V_\alpha$
 is hypercyclic on Fr\'{e}chet space
$C_0([0,1]):=\{u\in C([0,1]): u(0)=0\}$, endowed with the system of seminorms
$$
\|u\|_k=\max\limits_{t\in[0,1-1/(k+1)]}|u(t)|,\qquad k=1,2,\dots.
$$
If the operator $V_\alpha$ is defined on $L^p(0,1)$ ($1\le
p<\infty$) then  $\sigma(V_\alpha^*)$ is an infinite set and hence (see [3])
$V_\alpha$ cannot be even supercyclic on $L^p(0,1)$ .
 }

{\bf Acknowledgments.} I am grateful to Professor J. Zem\'{a}nek
for encouraging  me to study  the operator $V_{x^\alpha}$ and to
Professor M.M. Malamud for helpful remarks.

\references{14} {
 \item{[1]} G. E. Andrews, R. Askey  and R. Roy,
  {\it Laguerre Polynomials}\/, \S 6.2 in : Special
Functions, Cambridge University Press, 1999,
 282--293.

\item{[2]} A. Brown, P. R. Halmos and A. L. Shields, {\it Ces\'{a}ro
operators}, Acta Sci.Math. (Szeged) 26 (1965), 125--137.

\item{[3]} D. Herrero, {\it Limits of hypercyclic and supercyclic operators},
J. Funct. Anal. 99, No.1 (1991),   179--190.

\item{[4]} G. Herzog and A. Weber, {\it A class of hypercyclic
 Volterra composition operators}, Demonstratio Mathematica, v.XXXIX, No.2, (2006),
 465--468.

\item{[5]} Yu. I. Lyubich, {\it Composition of integration and
substitution}, in: Linear and Complex Analysis Problem Book,
Lecture Notes in Math. 1043, Springer,  Berlin, 1984, 249--250.

\item{[6]} Yu. I. Lyubich, {\it Linear Functional Analysis},
Moscow, 1988 (in Russian).

\item{[7]} C. M\H{u}ntz,
 {\it \H{U}ber den Approximationsatz von Weierstrass},
 H.A. Schwartz Feischrift, Berlin, 1914.

\item{[8]} B. Sz.-Nagy and C. Foias, {\it  Harmonic Analysis of
Operators on Hilbert Space}, Academiai Kiado, Budapest, 1970.

\item{[9]} G. P\'{o}lya and  G. Szeg\H{o}, {\it Problems and Theorems
in Analysis}, Springer, 1972.

\item{[10]} H. Radjavi and P. Rosenthal, {\it  Invariant Subspaces},
Springer, 1973.

\item{[11]} O. Sz\'{a}sz, {\it \H{U}ber die Approximation Steliger
Funktionen durch Lineare Aggregate von Potenzen}\/, Math.Ann. 77
(1916), 482--496.

\item{[12]} Yu Sun Tong, {\it Quasinilpotent integral operators},
Acta Math. Sinica 32 (1989), 727-735 (in Chinese).

\item{[13]}  R. Whitley, {\it The spectrum of a Volterra composition
operator}\/, Integral Equations Operator Theory 10 (1987),
146--149.

\item{[14]}  M. Zima, {\it A certain fixed point theorem and its
applications to integral-functional equations}\/, Bull. Austral.
Math. Soc. 46 (1992), 179--186. }
\end{document}